\documentclass[10pt]{article}

\usepackage{indentfirst}
\usepackage[a4paper,margin=1in]{geometry}
\usepackage{fontspec}
\usepackage{microtype}
\usepackage{amsmath,amssymb,amsthm,mathtools}
\usepackage{booktabs,array}
\usepackage{hyperref}
\hypersetup{colorlinks=true,linkcolor=blue,citecolor=blue,urlcolor=blue}
\usepackage{enumitem}
\setlist{nosep}
\allowdisplaybreaks

\newtheorem{theorem}{Theorem}[section]
\newtheorem{lemma}{Lemma}[section]
\newtheorem{proposition}{Proposition}[section]
\newtheorem{corollary}{Corollary}[section]

\theoremstyle{remark}
\newtheorem{remark}{Remark}[section]

\newcommand{\ex}{\operatorname{ex}}
\newcommand{\cC}{\mathcal{C}}

\title{Counting shorter even cycles in graphs forbidding an even cycle}
\author{Zhuoran HAN, Yaojun CHEN\footnote{Corresponding author. Email: yaojunc@nju.edu.cn}
\\
{\small School of Mathematics, Nanjing University, Nanjing 210093, China}
}
\date{}

\begin{document}
\maketitle

\begin{abstract}
For graphs $H$ and $F$, the generalized Tur\'an number $\ex(n,H,F)$ is the maximum
number of copies of $H$ in an $n$-vertex $F$-free graph. Let $C_m$ be a cycle on $m$ vertices.
Liu and Wu (2026) proved that, for fixed integers $k>\ell\ge2$,
\[
\ex(n,C_{2\ell},C_{2k})
=\left(\frac{(k-1)_\ell}{2\ell}+o(1)\right)n^\ell,
\]
where $(k-1)_\ell=(k-1)\cdots(k-\ell)$.
In this paper, we determine the value of $\ex(n,C_{2\ell},C_{2k})$
for all
sufficiently large $n$ and the unique extremal graph.

\end{abstract}

\medskip

\section{Introduction}
All graphs in this paper are finite, simple, and undirected, and all
subgraph copies are not necessarily induced.  For graphs $H$ and $G$,
let $N(H,G)$ denote the number of unlabeled copies of $H$ in $G$. The generalized Tur\'an number is
defined as:
\[
 \ex(n,H,F)=\max\{N(H,G): |V(G)|=n,\ G\text{ is }F\text{-free}\}.
\]
Write $K_t$, $I_t$, and $C_t$ for the complete graph,
the empty graph, and the cycle on $t$ vertices, respectively.
If $G_1$ and $G_2$ are vertex-disjoint, then $G_1\vee G_2$ denotes their
join.  For every nonnegative integer $r$, set
\[
 (x)_r=x(x-1)\cdots(x-r+1),\qquad (x)_0=1,
\]
and, when $x$ is a nonnegative integer with $x<r$, adopt the convention
$(x)_r=0$.

The systematic study of generalized Tur\'an numbers was initiated by
Alon and Shikhelman~\cite{AlonShikhelman}.  Cycle-counting problems in
this framework were further investigated by Gishboliner and
Shapira~\cite{GishbolinerShapira} and by Gerbner, Gy\H{o}ri, Methuku,
and Vizer~\cite{GerbnerGyoriMethukuVizer}. In particular, Gerbner, Gy\H{o}ri, Methuku, and Vizer \cite{GerbnerGyoriMethukuVizer} proved that
$\ex(n,C_{2\ell},C_{2k})=\Theta(n^\ell)$ for
$k>\ell\ge2$, and determined the leading term
when $\ell=2$. Liu and Wu~\cite{LiuWu} subsequently
determined the leading constant for $k>\ell\ge3$.
Their result, which also covers $\ell=2$, states that
\[
\ex(n,C_{2\ell},C_{2k})
=\left(\frac{(k-1)_\ell}{2\ell}+o(1)\right)n^\ell
\]
for every fixed pair $k>\ell\ge2$.

Only a few exact results were previously known when a single even cycle
is forbidden.  Gy\H{o}ri, He, Lv, Salia, Tompkins, Varga, and Zhu~\cite{GyoriEtAl}
proved that, for all sufficiently large $n$,
\[
 \ex(n,C_4,C_6)=\binom{n-2}{2}+2,
\]
with unique extremal graph
$K_2\vee(K_2\cup I_{n-4})$.  More recently, Chen and
Deng~\cite{ChenDeng} established
\[
 \ex(n,C_6,C_8)=6\binom{n-3}{3}+12(n-5),
\]
and showed that $K_3\vee(K_2\cup I_{n-5})$ is the unique extremal
graph.

A closely related exact problem concerns graphs of bounded
circumference.  Let $\cC_{\ge L}$ denote the family of all cycles of
length at least $L$. For $\cC_{\ge L}$-free graphs, Zhu, Gy\H{o}ri, He, Lv, Salia, and Xiao~\cite{ZhuEtAl} determined
the corresponding extremal number for $C_4$, and Zhao and
Wang~\cite{ZhaoWang} proved the analogue for every fixed even cycle $C_{2s}$ with $s\ge3$. However,
being $\cC_{\ge2k}$-free is strictly stronger than being merely
$C_{2k}$-free.  Thus those results do not rule out a graph containing
longer cycles but no $C_{2k}$ from having more copies of $C_{2\ell}$.

Set
\[
 a:=k-1,\qquad b:=n-a,
\]
and, for $n\ge k+1$, define
\begin{equation*}
 H_{n,k}:=K_a\vee\bigl(K_2\cup I_{b-2}\bigr).
\end{equation*}
In this paper, we try to count the maximum number of $C_{2\ell}$ in a $C_{2k}$-free graph for $k>\ell\ge2$.
Our main result gives the exact value and the unique extremal graph.

\begin{theorem}\label{thm:main}
For every fixed pair of integers $k>\ell\ge2$, there exists
$n_0=n_0(k,\ell)$ such that, for every $n\ge n_0$,
\[
 \ex(n,C_{2\ell},C_{2k})=N(C_{2\ell},H_{n,k}),
\]
where
\begin{align}
 N(C_{2\ell},H_{n,k})
 ={}&\sum_{j=0}^{\ell}
 \frac{(a)_{2\ell-j}}{2(2\ell-j)}
 \binom{2\ell-j}{j}(b)_j +
 \sum_{j=2}^{\ell}
 (a)_{2\ell-j}
 \binom{2\ell-j-1}{j-2}(b-2)_{j-2}.
 \label{eq:exact-formula}
\end{align}
Moreover, $H_{n,k}$ is the unique extremal graph.
\end{theorem}

When $k=\ell+1$, formula~\eqref{eq:exact-formula} simplifies as follows, which extends the results of Gy\H{o}ri, He, Lv, Salia, Tompkins, Varga, and Zhu \cite{GyoriEtAl}, and Chen and Deng~\cite{ChenDeng}.

\begin{corollary}\label{cor:adjacent}
For every fixed integer $\ell\ge2$ and all sufficiently large $n$,
\[
 \ex(n,C_{2\ell},C_{2\ell+2})
 =\frac{\ell!(\ell-1)!}{2}\binom{n-\ell}{\ell}
  +\ell!(\ell-1)!\binom{n-\ell-2}{\ell-2}.
\]
The unique extremal graph is $K_\ell\vee\bigl(K_2\cup I_{n-\ell-2}\bigr)$.
\end{corollary}
The proof of Theorem \ref{thm:main} is organized as follows.  Section~\ref{sec:candidate} verifies
that $H_{n,k}$ is $C_{2k}$-free and derives
\eqref{eq:exact-formula}. In Section~\ref{sec:counting-stability}, we first discuss a counting tool and then apply it to show that for any near-extremal graph $G$, there is a set $A\subseteq V(G)$ of size $a$ satisfying $\left|\bigcap_{v\in A}N_G(v)\right|=n-o(n)$. Section~\ref{sec:exceptional} compares the near-extremal graph with $H_{n,k}$, rules out the existence of exceptional vertices, and thereby completes the proof of the main theorem.

Throughout the paper, $k$ and $\ell$ are fixed, and the implicit
constants in all $O(\cdot)$ and $o(\cdot)$ terms depend only on
$k$ and $\ell$. We write $N_G(v)$ and $d_G(v)$ for the neighborhood and
degree of a vertex $v$, and $e(G)$ for the number of edges of $G$.  For
$X\subseteq V(G)$, the subgraph induced by $X$ is denoted by $G[X]$.

\section{The candidate graph}\label{sec:candidate}

\begin{proposition}\label{prop:candidate}
Let $n\ge k+1$. Then $H_{n,k}$ is $C_{2k}$-free, and its number of
copies of $C_{2\ell}$ is given by \eqref{eq:exact-formula}.
\end{proposition}

\begin{proof}
Let $A$ be the vertex set of the copy of $K_a$ in $H_{n,k}$, and let
$B=V(H_{n,k})\setminus A$.  Thus $|A|=a$, $|B|=b$, all edges between
$A$ and $B$ are present, and $H_{n,k}[B]$ consists of a single edge.

Consider a cycle $C$ in $H_{n,k}$.  Suppose that $C$ uses $s$ vertices
of $A$ and $j$ vertices of $B$, and let $e_A$, $e_B$, and $e_\times$
be the numbers of edges of $C$ lying inside $A$, inside $B$, and between
the two parts, respectively.  Summing degrees on the two sides of the
cycle gives
\[
 2s=2e_A+e_\times,
 \qquad
 2j=2e_B+e_\times,
\]
and hence
\[
 j-s=e_B-e_A.
\]
If $C$ had length $2k$, then $s\le a=k-1$, so
$j-s=2k-2s\ge2$.  On the other hand, $e_B\le1$ and $e_A\ge0$, giving
$j-s\le1$, a contradiction.  Thus $H_{n,k}$ is $C_{2k}$-free.

Now let $C$ be a copy of $C_{2\ell}$.  Since $s+j=2\ell$ and
$j-s=e_B-e_A\le1$, we have $0\le j\le\ell$ and
$s=2\ell-j$.

First count the cycles that do not use the unique edge in $B$.  Record,
in cyclic order and with a distinguished starting point, the vertices
of $A$ encountered along an orientation of the cycle.  There are
$(a)_s$ choices.  Choose $j$ of the $s$ cyclic gaps and place distinct
vertices of $B$ in these gaps; this gives
$\binom{s}{j}(b)_j$ choices.  Every unoriented cycle has exactly $2s$
such representations.  Hence the number of these cycles is
\[
 \frac{(a)_s}{2s}\binom{s}{j}(b)_j.
\]
Substituting $s=2\ell-j$ and summing over $0\le j\le\ell$ gives the
first sum in \eqref{eq:exact-formula}.

Next count the cycles that use the unique edge $xy$ in $B$.  Regard
$xy$ as an oriented two-vertex block and the remaining $j-2$ vertices
of $B$ as singleton blocks.  The resulting $j-1$ external blocks must
occupy distinct gaps of the cyclic core sequence.  After choosing and
orienting the two-vertex block and placing the remaining distinct
external vertices, and then dividing by the $2s$ representations of
each unoriented cycle, the number obtained is
\[
 (a)_s\binom{s-1}{j-2}(b-2)_{j-2}.
\]
Substituting $s=2\ell-j$ and summing over $2\le j\le\ell$ gives the
second sum in \eqref{eq:exact-formula}.
\end{proof}

The term with $j=\ell$ in the first sum in \eqref{eq:exact-formula} has degree $\ell$ in $n$, whereas all remaining terms have degree at most $\ell-1$.
Consequently,
\begin{equation}\label{eq:candidate-asymptotic}
	N(C_{2\ell},H_{n,k})
	=\frac{(a)_\ell}{2\ell}(b)_\ell
	+O_{k,\ell}(n^{\ell-1})
	=\frac{(a)_\ell}{2\ell}n^\ell
	+O_{k,\ell}(n^{\ell-1}).
\end{equation}

\section{Counting tools and near-extremal structure}
\label{sec:counting-stability}

In this section, we adapt the high-codegree truncation and
neighborhood-clique construction of Liu and Wu~\cite{LiuWu}
to bound the number of target cycles using a small auxiliary
graph. We then characterize the equality case of the
compatible-ordering bound and deduce that every
asymptotically extremal graph contains a set of $k-1$
vertices with $n-o(n)$ common neighbors.

\subsection{Counting tools}\label{subsec:counting-tools}

Let $G$ be an $n$-vertex $C_{2k}$-free graph.  For distinct vertices
$u,v\in V(G)$, define their codegree by
\[
 f(u,v):=|N_G(u)\cap N_G(v)|.
\]
For a real number $q>0$, let $H_q$ be the graph with edge set
\[
 E(H_q)=\{uv:f(u,v)\ge q\},
\]
with isolated vertices deleted.  Thus $uv\in E(H_q)$ records only a
large codegree in $G$; it does not assert that $uv\in E(G)$.

We begin with the classical even-cycle extremal estimate.

\begin{theorem}[Bondy and Simonovits~\cite{BondySimonovits}]
\label{thm:Bondy-Simonovits}
For every integer $k\ge2$, every $N$-vertex $C_{2k}$-free graph $J$
satisfies
\[
 e(J)=O_k(N^{1+1/k}).
\]
\end{theorem}

Then we recall the following high-codegree estimates of Liu and Wu.

\begin{lemma}[Liu and Wu~\cite{LiuWu}]
\label{lem:high-codegree}
Let $k\ge3$, and let $G$ be an $n$-vertex $C_{2k}$-free graph.  For every
real $q\ge2k$, the following statements hold.
\begin{enumerate}[label=\textup{(\roman*)}]
\item $H_q$ is $C_k$-free.
\item For every subgraph $J\subseteq H_q$,
\[
 \sum_{uv\in E(J)} f(u,v)
 \le \frac{k-2}{2}\sum_{v\in V(J)}d_G(v)
 =O_k\!\left(n|V(J)|^{1/k}\right).
\]
\item
\[
 e(H_q)=O_k\!\left((n/q)^{k/(k-1)}\right).
\]
\item For every real $p$ with $2k\le p\le n$,
\[
 \sum_{\substack{uv\in E(H_{2k})\\ f(u,v)\le p}}
 f(u,v)^2
 =O_k\!\left(n^{k/(k-1)}p^{(k-2)/(k-1)}\right).
\]
In particular,
\[
 \sum_{uv\in E(H_{2k})}f(u,v)^2=O_k(n^2).
\]
\end{enumerate}
\end{lemma}

We will also use the following cyclic form of the Cauchy--Schwarz
inequality.

\begin{lemma}
	\label{lem:cyclic-CS}
	Let $X$ be a finite set, let $r\ge2$, and let
	$a_1,\ldots,a_r:X\times X\to[0,\infty)$.  With $x_{r+1}=x_1$,
	\[
	\sum_{x_1,\ldots,x_r\in X}\prod_{i=1}^r a_i(x_i,x_{i+1})
	\le
	\prod_{i=1}^r
	\left(\sum_{x,y\in X}a_i(x,y)^2\right)^{1/2}.
	\]
\end{lemma}

\begin{proof}
	We argue by induction on $r$.  The case $r=2$ is the ordinary
	Cauchy--Schwarz inequality.  Suppose $r\ge3$ and define
	\[
	b(x,z):=\sum_{y\in X}a_1(x,y)a_2(y,z).
	\]
	For fixed $x,z$, Cauchy--Schwarz gives
	\[
	b(x,z)^2
	\le
	\left(\sum_{y\in X}a_1(x,y)^2\right)
	\left(\sum_{y\in X}a_2(y,z)^2\right).
	\]
	Summing over $x,z$ yields
	\[
	\sum_{x,z\in X}b(x,z)^2
	\le
	\left(\sum_{x,y\in X}a_1(x,y)^2\right)
	\left(\sum_{y,z\in X}a_2(y,z)^2\right).
	\]
	The left-hand side in the statement of the lemma equals
	\[
	\sum_{x_1,x_3,\ldots,x_r}
	b(x_1,x_3)a_3(x_3,x_4)\cdots a_r(x_r,x_1).
	\]
	Applying the induction hypothesis to $b,a_3,\ldots,a_r$ and then using
	the preceding bound proves the claim.
\end{proof}

Write a copy of $C_{2\ell}$ as
\[
 v_0x_0v_1x_1\cdots v_{\ell-1}x_{\ell-1}v_0,
\]
with indices taken modulo $\ell$.  The two sequences
$v_0,\ldots,v_{\ell-1}$ and $x_0,\ldots,x_{\ell-1}$ are called the two
\emph{alternating sides} of the cycle.  An alternating side is called
\emph{high-codegree} if every pair of consecutive vertices on that side
has codegree at least $2k$ in $G$.  A target cycle is called \emph{bad}
if neither alternating side is high-codegree.

\begin{lemma}\label{lem:bad-cycles}
Let $k>\ell\ge2$, and let $G$ be an $n$-vertex $C_{2k}$-free graph.  The
number of bad copies of $C_{2\ell}$ in $G$ is $o(n^\ell)$.
\end{lemma}

\begin{proof}
For $\ell\ge3$, the assertion is exactly the exceptional-cycle estimate
of Liu and Wu~\cite{LiuWu}.  Their proof uses a path of length
$2\ell-5$, so the case $\ell=2$ requires a separate argument.

Assume $\ell=2$.  A bad four-cycle has both pairs of opposite vertices
of codegree less than $2k$.  By Theorem~\ref{thm:Bondy-Simonovits}, every
$s$-vertex induced subgraph of $G$ has average degree
$O_k(s^{1/k})\le O_k(n^{1/k})$.  Repeatedly deleting a minimum-degree
vertex therefore yields an ordering of $V(G)$ in which every vertex has
at most $D=O_k(n^{1/k})$ neighbors appearing later in the order.

Fix a bad four-cycle $x-y-z-w-x$ and let $x$ be its first vertex in this
ordering.  There are at most $n$ choices for $x$, at most
$\binom{D}{2}$ choices for the unordered pair $\{y,w\}$, and fewer than
$2k$ choices for $z$, because $f(y,w)<2k$.  Ignoring the remaining
distinctness and adjacency conditions only enlarges the count.  Hence
the number of bad four-cycles is at most
\[
 n\binom{D}{2}(2k-1)=O_k(n^{1+2/k})=o(n^2),
\]
since $k\ge3$.
\end{proof}

Set
\[
 q:=\frac{n}{\log n},
 \qquad
 W:=V(H_q),
 \qquad
 F:=H_{2k}[W].
\]
For all sufficiently large $n$, we have $q\ge2k$.  By
Lemma~\ref{lem:high-codegree}\textup{(i),(iii)}, the graph $F$ is
$C_k$-free and
\[
 |W|\le2e(H_q)=O_k\!\left((\log n)^{k/(k-1)}\right).
\]

\begin{lemma}
\label{lem:localization}
Let $k>\ell\ge2$, and let $G$ be an $n$-vertex $C_{2k}$-free graph. All but $o(n^\ell)$ copies of $C_{2\ell}$ in $G$ have a high-codegree alternating side entirely contained in $W$.
\end{lemma}

\begin{proof}
By Lemma~\ref{lem:bad-cycles}, it remains to count cycles having a
high-codegree alternating side for which a chosen such side is not
entirely contained in $W$.  Define symmetric weight functions on
$V(G)\times V(G)$ by
\[
 \omega(u,v)=
 \begin{cases}
  f(u,v),&uv\in E(H_{2k}),\\
  0,&\text{otherwise},
 \end{cases}
\]
and
\[
 \omega_{<q}(u,v)=
 \begin{cases}
  \omega(u,v),&0<\omega(u,v)<q,\\
  0,&\text{otherwise}.
 \end{cases}
\]
Lemma~\ref{lem:high-codegree}\textup{(iv)} gives
\[
 \sum_{u,v}\omega(u,v)^2=O_k(n^2)
\]
and
\[
 \sum_{u,v}\omega_{<q}(u,v)^2
 =O_k\!\left(\frac{n^{2}}{(\log n)^{(k-2)/(k-1)}}\right)
 =o(n^2).
\]

If the chosen high-codegree alternating side is not entirely
contained in $W$, then it contains a vertex $v\notin W$.
Since $W=V(H_q)$ and $H_q$ has no isolated vertices,
$v$ is not incident to any edge of $H_q$.
Every consecutive pair on the chosen side has codegree at least
$2k$, while any pair containing $v$ has codegree less than $q$.
Hence, after a cyclic relabeling, the cycle can be written as
\[
v_1x_1v_2x_2\cdots v_\ell x_\ell v_1
\]
with $2k\le f(v_1,v_2)<q$, and consequently $\omega_{<q}(v_1,v_2)>0$.
  Once
$v_1,\ldots,v_\ell$ are fixed, each $x_i$ must be a common neighbor of
$v_i$ and $v_{i+1}$, so the number of choices is at most
\[
 \omega_{<q}(v_1,v_2)\omega(v_2,v_3)\cdots\omega(v_\ell,v_1).
\]

Ignoring all distinctness conditions, we sum over all
ordered tuples $(v_1,\ldots,v_\ell)\in V(G)^\ell$.
This summation already accounts for every possible choice
of the high-codegree alternating side and its cyclic labeling.
By Lemma~\ref{lem:cyclic-CS}, the total number of
such cycles is at most
\[
\begin{aligned}
&\sum_{v_1,\ldots,v_\ell\in V(G)}
\omega_{<q}(v_1,v_2)
\omega(v_2,v_3)\cdots\omega(v_\ell,v_1)\\
&\qquad\le
\left(\sum_{u,v}\omega_{<q}(u,v)^2\right)^{1/2}
\left(\sum_{u,v}\omega(u,v)^2\right)^{(\ell-1)/2}\\
&\qquad=o(n^\ell),
\end{aligned}
\]
where the last equality follows from the preceding estimates.

\end{proof}

For each $x\in V(G)$, let
\[
 P_x:=N_G(x)\cap W.
\]
We regard the empty set and singleton sets as cliques, and define
\[
 Q_x:=
 \begin{cases}
  P_x,&F[P_x]\text{ is complete},\\
  \varnothing,&\text{otherwise}.
 \end{cases}
\]
Then $Q_x$ is a clique of $F$ and $Q_x\subseteq N_G(x)$.  Different
labels $x$ may determine the same clique.

For a label sequence
$\mathbf{x}=(x_1,\ldots,x_\ell)\in V(G)^\ell$, consider ordered
$\ell$-tuples of distinct vertices
$(v_1,\ldots,v_\ell)\in W^\ell$.  Such a tuple is said to be
\emph{compatible} with $\mathbf{x}$ if
\[
 \{v_i,v_{i+1}\}\subseteq Q_{x_i}
 \quad (1\le i\le\ell),
 \qquad v_{\ell+1}=v_1.
\]
Let $\kappa(\mathbf{x})$ be the number of compatible ordered tuples. For $\ell\ge3$, every compatible tuple $(v_1,\ldots,v_\ell)$ forms a cycle in $F$.
When $\ell=2$, compatibility means $\{v_1,v_2\}\subseteq Q_{x_1}\cap Q_{x_2}$.
We do not require the labels $x_1,\ldots,x_\ell$ to be
distinct from one another or from the $v_i$, since
allowing such repetitions can only increase the count.

\begin{lemma}\label{lem:vertex-labeled}
With the preceding notation,
\begin{equation}\label{eq:vertex-labeled-count}
 2\ell N(C_{2\ell},G)
 \le
 \sum_{\mathbf{x}\in V(G)^\ell}\kappa(\mathbf{x})+o(n^\ell).
\end{equation}
\end{lemma}

\begin{proof}
Let
\[
 D:=\{x\in V(G):F[P_x]\text{ is not complete}\}.
\]
For each $x\in D$, there exist distinct $u,v\in P_x$
such that $uv\notin E(F)$. Since $u,v\in W$ and
$F=H_{2k}[W]$, we have $f(u,v)<2k$. Moreover, $x$ is
a common neighbor of $u$ and $v$. Hence
\[
 |D|
 \le \sum_{\substack{\{u,v\}\subseteq W\\uv\notin E(F)}}
 f(u,v)
 \le (2k-1)\binom{|W|}{2}
 =O_k\!\left((\log n)^{2k/(k-1)}\right).
\]

By Lemma~\ref{lem:localization}, all but $o(n^\ell)$
copies of $C_{2\ell}$ have a high-codegree alternating
side contained in $W$. For each such cycle, choose
one such side and fix a cyclic labeling
\[
 C=v_1x_1v_2x_2\cdots v_\ell x_\ell v_1,
 \qquad v_1,\ldots,v_\ell\in W.
\]
We further discard cycles for which some $x_i$ belongs
to $D$. To bound their number, fix a position $i$.
There are at most $|W|^\ell$ choices for
$(v_1,\ldots,v_\ell)$, at most $|D|$ choices for $x_i$,
and at most $n^{\ell-1}$ choices for the remaining
$x_j$. Thus, ignoring all distinctness and adjacency
conditions, the number of additional discarded cycles
is at most
\[
\begin{aligned}
 \ell |W|^\ell |D|n^{\ell-1}
 &=O_{k,\ell}\!\left(
 n^{\ell-1}(\log n)^{(\ell+2)k/(k-1)}
 \right)=o(n^\ell).
\end{aligned}
\]

For every remaining cycle, we have
$v_i,v_{i+1}\in P_{x_i}$ for each $i$, where
$v_{\ell+1}=v_1$. Since $x_i\notin D$,
we have $Q_{x_i}=P_{x_i}$, and hence
\[
 \{v_i,v_{i+1}\}\subseteq Q_{x_i}
 \qquad (1\le i\le\ell).
\]
Therefore $(v_1,\ldots,v_\ell)$ is compatible
with $(x_1,\ldots,x_\ell)$.

With the chosen alternating side fixed, each remaining
unoriented cycle has exactly $2\ell$ distinct compatible
representations, obtained by choosing one of the $\ell$
starting vertices on that side and one of the two
directions. Each representation consists of the two
ordered tuples $(x_1,\ldots,x_\ell)$ and
$(v_1,\ldots,v_\ell)$, which together uniquely determine
the original cycle. Thus representations arising from
different cycles are distinct.

Since the total number of discarded cycles is
$o(n^\ell)$, we obtain
\[
 2\ell N(C_{2\ell},G)
 \le
 \sum_{\mathbf{x}\in V(G)^\ell}\kappa(\mathbf{x})
 +o(n^\ell),
\]
as required.
\end{proof}

\subsection{The structure of near-extremal graphs}
\label{subsec:stability}

The vertex-labeled count reduces the global cycle problem to the number
of compatible orderings associated with each label sequence.  The next
lemma gives the optimal bound for one sequence and characterizes its
equality case.

\begin{lemma}
\label{lem:rigidity}
Let $a=k-1$.  For every label sequence
$\mathbf{x}=(x_1,\ldots,x_\ell)$,
\[
 \kappa(\mathbf{x})\le(a)_\ell.
\]
Equality holds if and only if there exists an $a$-vertex clique $X$ in
$F$ such that
\[
 Q_{x_1}=\cdots=Q_{x_\ell}=X.
\]
\end{lemma}

\begin{proof}
If no compatible tuple exists, the inequality is immediate.  Otherwise,
fix a compatible tuple $(v_1,\ldots,v_\ell)$ and set
\[
 Z:=\bigcup_{i=1}^{\ell}Q_{x_i}.
\]
We first show that $|Z|\le k-1$.  Suppose instead that $|Z|\ge k$.
Since the $v_i$ are distinct, we may choose $k-\ell$ distinct vertices
from $Z\setminus\{v_1,\ldots,v_\ell\}$.  Assign each chosen vertex to
an index $i$ for which it belongs to $Q_{x_i}$.  Because $Q_{x_i}$ is a
clique containing $v_i$ and $v_{i+1}$, the vertices assigned to $i$ can
be inserted, in any order, between $v_i$ and $v_{i+1}$.  The resulting
paths are internally disjoint and concatenate to a cycle of length $k$
in $F$, contradicting the fact that $F$ is $C_k$-free.  When $\ell=2$,
the construction gives two internally disjoint paths joining $v_1$ and
$v_2$; at least one path has an internal vertex, so their union is again
a simple $C_k$.

Consequently $|Z|\le a$, and every compatible tuple is an ordered
$\ell$-tuple of distinct vertices of $Z$.  Therefore
\[
 \kappa(\mathbf{x})\le(|Z|)_\ell\le(a)_\ell.
\]
If equality holds, then $|Z|=a$ and every ordered $\ell$-tuple of
distinct vertices of $Z$ is compatible.  Fix $i$ and $z\in Z$.  Since
$a\ge\ell$, there exists such a tuple with $z$ in position $i$.
Compatibility then gives $z\in Q_{x_i}$.  Thus $Z\subseteq Q_{x_i}$,
and the reverse inclusion is automatic from the definition of $Z$.
Hence $Q_{x_i}=Z$ for every $i$.  Conversely, if all $Q_{x_i}$ are the
same $a$-vertex clique, then all $(a)_\ell$ ordered $\ell$-tuples of
distinct vertices of that clique are compatible.
\end{proof}

Combining Lemmas~\ref{lem:vertex-labeled} and~\ref{lem:rigidity} gives
\[
 N(C_{2\ell},G)
 \le\left(\frac{(a)_\ell}{2\ell}+o(1)\right)n^\ell,
\]
which also recovers the asymptotic upper bound of Liu and Wu.  We will use the following uniform coarse
consequence in Section~\ref{sec:exceptional}.

\begin{corollary}
\label{cor:coarse-bound}
There exists a constant $C_0=C_0(k,\ell)>0$ such that every
$N$-vertex $C_{2k}$-free graph $J$ satisfies
\[
 N(C_{2\ell},J)\le C_0N^\ell.
\]
\end{corollary}

\begin{proof}
For all sufficiently large $N$, the preceding asymptotic upper bound
holds with one fixed constant.  Enlarging that constant to dominate
$N(C_{2\ell},K_N)/N^\ell$ for the finitely many remaining values of $N$
proves the assertion.
\end{proof}

\begin{lemma}
\label{lem:common-neighborhood}
Fix $k>\ell\ge2$ and let $a=k-1$.  Suppose that $G=G_n$ is a sequence
of $n$-vertex $C_{2k}$-free graphs satisfying
\[
 N(C_{2\ell},G)
 \ge\left(\frac{(a)_\ell}{2\ell}-o(1)\right)n^\ell.
\]
Then there exists a set $A\subseteq V(G)$ of size $a$ such that
\[
 \left|\bigcap_{v\in A}N_G(v)\right|=n-o(n).
\]
\end{lemma}

\begin{proof}
Construct $F$ and the family $(Q_x)_{x\in V(G)}$ as in
Subsection~\ref{subsec:counting-tools}, and write
\[
 L:=(a)_\ell.
\]
By \eqref{eq:vertex-labeled-count} and the hypothesis,
\[
 \sum_{\mathbf{x}\in V(G)^\ell}\kappa(\mathbf{x})
 \ge Ln^\ell-o(n^\ell).
\]
For each $a$-vertex clique $X$ of $F$, let
\[
 m_X:=|\{x\in V(G):Q_x=X\}|,
\]
and let $M:=\max_Xm_X$, with $M=0$ if $F$ contains no $a$-vertex
clique.  Let $T$ be the number of label sequences $\mathbf{x}$ for which
$\kappa(\mathbf{x})=L$.  By Lemma~\ref{lem:rigidity}, all labels in such
a sequence correspond to one common $a$-vertex clique.  After choosing
$x_1$ arbitrarily, there are at most $M$ choices for $x_2$, and at most
$n^{\ell-2}$ choices for the remaining labels.  Hence
\[
 T\le Mn^{\ell-1}.
\]
Since every $\kappa(\mathbf{x})$ is a nonnegative integer at most $L$,
all sequences not counted by $T$ contribute at most $L-1$.  Therefore
\begin{align*}
 \sum_{\mathbf{x}\in V(G)^\ell}\kappa(\mathbf{x})
 &\le LT+(L-1)(n^\ell-T)\\
 &=(L-1)n^\ell+T\\
 &\le(L-1)n^\ell+Mn^{\ell-1}.
\end{align*}
Comparing this estimate with the preceding lower bound gives
\[
 (n-M)n^{\ell-1}=o(n^\ell),
\]
so $M=n-o(n)$.

Choose an $a$-vertex clique $A$ of $F$ attaining the maximum
multiplicity $M$.  For every $x$ with $Q_x=A$, we have
$A\subseteq N_G(x)$.  Thus $x$ is a common neighbor of all vertices in $A$, and
\[
 \left|\bigcap_{v\in A}N_G(v)\right|
 \ge|\{x:Q_x=A\}|=M=n-o(n),
\]
as required.
\end{proof}

\begin{remark}
 The set $A$ induces a clique in $F$, but is not asserted at this stage to induce a clique in $G$.
\end{remark}

\section{Eliminating exceptional vertices and proving the main theorem}
\label{sec:exceptional}

Let $G$ be a $C_{2k}$-free graph satisfying the hypothesis of
Lemma~\ref{lem:common-neighborhood}, and let $A$ be the $a$-vertex set
given by that lemma. Define
\[
B:=\bigcap_{v\in A}N_G(v),
\qquad
R:=V(G)\setminus(A\cup B),
\qquad
r:=|R|=o(n).
\]
Since $G$ has no loops, $A\cap B=\varnothing$. All edges between $A$
and $B$ are present, whereas every vertex of $R$ misses at least one
edge to $A$.

\begin{lemma}\label{lem:local-structure}
	For all sufficiently large $n$,
	\[
	e(G[B])\le1,
	\qquad
	d_B(x)\le1\quad\text{for every }x\in R.
	\]
	Moreover, there exists $U\subseteq B$ such that
	\[
	N_G(u)=A\quad\text{for every }u\in U
	\]
	and
	\[
	|V(G)\setminus U|\le2r+a+2.
	\]
\end{lemma}

\begin{proof}
	Suppose first that $G[B]$ contains two adjacent edges. Regard the
	resulting three-vertex path as one path block and choose $a-1$ further
	vertices of $B$ as singleton blocks. If $G[B]$ contains two disjoint
	edges, regard them as two two-vertex path blocks and choose $a-2$
	further singleton blocks. Since $|B|=n-o(n)$, all required vertices
	exist for sufficiently large $n$. In either case there are $a$
	vertex-disjoint path blocks in $B$ containing altogether $a+2$
	vertices. Arrange the blocks cyclically and insert the $a$ vertices of
	$A$, one between each pair of consecutive blocks. All required joining
	edges are present, so this produces a cycle of length
	\[
	a+(a+2)=2k,
	\]
	a contradiction. Thus $e(G[B])\le1$.
	
	Now suppose that some $x\in R$ has two distinct neighbors
	$b_0,b_a\in B$. Write $A=\{a_1,\ldots,a_a\}$ and choose distinct
	vertices $b_1,\ldots,b_{a-1}$ from $B\setminus\{b_0,b_a\}$. Then
	\[
	xb_0a_1b_1a_2b_2\cdots a_ab_ax
	\]
	is a copy of $C_{2k}$, again a contradiction. Hence $d_B(x)\le1$ for
	every $x\in R$.
	
	Delete from $B$ every vertex adjacent to some vertex of $R$, and also
	delete the two endpoints of the possible unique edge in $G[B]$. Let
	$U$ be the remaining set. The first deletion removes at most $r$
	vertices and the second at most two. Every $u\in U$ is adjacent to all
	vertices of $A$, to no vertex of $R$, and to no other vertex of $B$.
	Thus $N_G(u)=A$, and
	\[
	|V(G)\setminus U|
	\le a+r+(r+2)=2r+a+2.
	\]
\end{proof}

Set
\[
m:=|U|,
\qquad
T:=V(G)\setminus U,
\qquad
t:=|T|.
\]
Then
\[
t\le2r+a+2=o(n),
\qquad
m=n-t\sim n.
\]
If $G[B]$ contains an edge $xy$, retain that edge. Otherwise, choose any
two distinct vertices $x,y\in B$. On the same vertex set as $G$, define
$G^\star$ by making $A$ complete, adding all edges between $A$ and
$V(G)\setminus A$, and retaining inside $V(G)\setminus A$ only the edge
$xy$. Then
\[
G^\star\cong H_{n,k},
\]
and Proposition~\ref{prop:candidate} shows that $G^\star$ is
$C_{2k}$-free.

\begin{lemma}\label{lem:gain-loss}
	If $r\ge1$, then, for all sufficiently large $n$,
	\[
	N(C_{2\ell},G)<N(C_{2\ell},G^\star).
	\]
\end{lemma}

\begin{proof}
	For a graph $J$, let $\mathfrak C(J)$ be the set of its copies of
	$C_{2\ell}$, and write
	\[
	\mathcal N:=\mathfrak C(G^\star)\setminus\mathfrak C(G),
	\qquad
	\mathcal L:=\mathfrak C(G)\setminus\mathfrak C(G^\star).
	\]
	We first count new cycles. For each $v\in R$, choose a vertex
	$a_1(v)\in A$ such that $va_1(v)\notin E(G)$, and fix distinct vertices
	\[
	a_2(v),\ldots,a_\ell(v)\in A\setminus\{a_1(v)\}.
	\]
	This is possible because $a=k-1\ge\ell$. For every ordered
	$(\ell-1)$-tuple of distinct vertices
	$(u_1,\ldots,u_{\ell-1})\in U^{\ell-1}$, the graph $G^\star$ contains
	the cycle
	\[
	va_1(v)u_1a_2(v)u_2\cdots u_{\ell-1}a_\ell(v)v.
	\]
	It is not a cycle of $G$, because $va_1(v)\notin E(G)$. Each such cycle
	has $v$ as its unique vertex of $R$. Starting at $v$ and following the
	fixed edge $va_1(v)$ recovers the ordered $U$-tuple, so all the cycles
	constructed in this way are distinct. Consequently,
	\[
	|\mathcal N|\ge r(m)_{\ell-1}.
	\]
	
	We next estimate the lost cycles. Every $C\in\mathcal L$ must contain
	an edge with both endpoints in $V(G)\setminus A$ other than $xy$, since
	$G^\star$ contains all edges inside $A$ and all edges between $A$ and
	its complement. Neither endpoint of this edge lies in $U$, because
	$N_G(u)=A$ for every $u\in U$.
	
	Let $j=|V(C)\cap U|$. If $j\ge1$, deleting these $j$ vertices from the
	cycle leaves $j$ nonempty paths. Each such path has both endpoints in
	$A$, since every vertex of $U$ has neighborhood exactly $A$. The path
	containing an edge absent from $G^\star$ has two distinct endpoints in
	$A$, while both endpoints of that absent edge lie outside $A$.
	Therefore this path contains at least four distinct vertices. Each of
	the other $j-1$ paths contains at least one vertex. Hence
	\[
	2\ell-j\ge4+(j-1)=j+3,
	\]
	which implies $j\le\ell-2$. The same conclusion is immediate when
	$j=0$. Thus every lost cycle contains at most $\ell-2$ vertices of $U$.
	
	For sufficiently large $n$, we have $m\ge t$. Partition $U$ into
	\[
	h:=\left\lceil\frac{m}{t}\right\rceil
	\]
	blocks, each containing at most $t$ vertices. Every lost cycle meets at
	most $\ell-2$ blocks. For any choice of $s\le\ell-2$ blocks, the induced
	graph on their union with $T$ has at most $(s+1)t$ vertices and remains
	$C_{2k}$-free. Therefore Corollary~\ref{cor:coarse-bound} and a union
	bound give
	\[
	|\mathcal L|
	\le \sum_{s=0}^{\ell-2}\binom hs C_0((s+1)t)^\ell
	=O_{k,\ell}(h^{\ell-2}t^\ell).
	\]
	Since $m\ge t$, we have $h\le2m/t$, and hence
	\[
	|\mathcal L|=O_{k,\ell}(t^2m^{\ell-2}).
	\]

	Combining the two preceding estimates yields
	\begin{equation}\label{eq:gain-loss-difference}
		N(C_{2\ell},G^\star)-N(C_{2\ell},G)
		\ge r(m)_{\ell-1}-O_{k,\ell}(t^2m^{\ell-2}).
	\end{equation}
	If $r\ge1$, then the bound on $t$ above gives
	\[
	t\le2r+a+2\le(a+4)r.
	\]
	Since $r=o(n)$ and $m\sim n$,
	\[
	t^2m^{\ell-2}=O_k(r^2m^{\ell-2})=o(rm^{\ell-1}),
	\]
	whereas
	\[
	r(m)_{\ell-1}=\Theta_\ell(rm^{\ell-1}).
	\]
	The right-hand side of \eqref{eq:gain-loss-difference} is therefore
	positive for all sufficiently large $n$.
\end{proof}

\begin{proof}[Proof of Theorem~\ref{thm:main}]
	Proposition~\ref{prop:candidate} shows that $H_{n,k}$ is admissible and
	gives its exact number of target cycles. Suppose that for infinitely
	many $n$ there exists an extremal $C_{2k}$-free graph $G_n$ not
	isomorphic to $H_{n,k}$. Consider such a sequence and write $G=G_n$.
	
	By extremality and \eqref{eq:candidate-asymptotic},
	\[
	N(C_{2\ell},G)
	\ge\left(\frac{(a)_\ell}{2\ell}-o(1)\right)n^\ell.
	\]
	Lemma~\ref{lem:common-neighborhood} therefore gives an $a$-vertex set
	$A$ whose common neighborhood $B$ has size $n-o(n)$. Use the notation
	introduced at the beginning of this section. If $r=|R|\ge1$ along an
	infinite subsequence, then Lemma~\ref{lem:gain-loss} produces a
	$C_{2k}$-free graph $G^\star\cong H_{n,k}$ with
	\[
	N(C_{2\ell},G^\star)>N(C_{2\ell},G),
	\]
	contradicting extremality. Thus $r=0$ for all sufficiently large $n$ in
	the sequence.
	
	Thus $V(G)=A\mathbin{\dot\cup}B$, all edges between $A$ and $B$ are
	present, and Lemma~\ref{lem:local-structure} gives $e(G[B])\le1$. Make
	$A$ complete, and, if $G[B]$ is empty, add an arbitrary edge $xy$ in
	$B$; otherwise, let $xy$ be its unique edge. Denote the resulting graph
	again by $G^\star$. Then
	\[
	G\subseteq G^\star,
	\qquad
	G^\star\cong H_{n,k}.
	\]
	
	Suppose $G\ne G^\star$. If $xy\notin E(G)$, choose arbitrary distinct
	vertices $u,v\in A$. Otherwise, $xy\in E(G)$ and some edge
	$uv\in\binom{A}{2}$ is missing from $G$; choose such an edge. Select
	distinct vertices
	\[
	a_3,\ldots,a_\ell\in A\setminus\{u,v\}
	\]
	and
	\[
	z_1,\ldots,z_{\ell-2}\in B\setminus\{x,y\},
	\]
	with both lists understood to be empty when $\ell=2$. The graph
	$G^\star$ contains the cycle
	\[
	uvz_1a_3z_2a_4\cdots z_{\ell-2}a_\ell xyu.
	\]
	For $\ell=2$, this is simply $uvxyu$. The cycle uses both $uv$ and
	$xy$, and therefore is absent from $G$ in either case. Since
	$G\subseteq G^\star$, it follows that
	\[
	N(C_{2\ell},G)<N(C_{2\ell},G^\star),
	\]
	again contradicting extremality. Consequently
	$G=G^\star\cong H_{n,k}$, completing the proof.
\end{proof}

\section{Concluding remarks}

In this paper, we have determined the exact value of
$\ex(n,C_{2\ell},C_{2k})$ for all fixed integers
$k>\ell\ge2$ and sufficiently large $n$, and proved
that $H_{n,k}$ is the unique extremal graph.

Several related questions remain of interest.
First, our result holds for sufficiently large $n$,
but the proof does not determine the smallest
possible threshold. This leads to the following
question.

\medskip
\noindent
\textbf{Question 1.}
For fixed integers $k>\ell\ge2$, what is the
smallest integer $n_0(k,\ell)$ such that
$H_{n,k}$ is the unique extremal graph for every
$n\ge n_0(k,\ell)$?

\medskip
Clearly, $n_0(k,\ell)\ge 2k$, since $K_{2k-1}$ is
$C_{2k}$-free and contains more copies of $C_{2\ell}$
than $H_{2k-1,k}$. It would be interesting to determine
whether the threshold can be expressed explicitly
in terms of $k$ and $\ell$.

Then, the common-neighborhood stability
established in this paper suggests a more precise
structural problem.

\medskip
\noindent
\textbf{Question 2.}
Can one obtain a sharp quantitative stability
theorem relating the deficit
\[
N(C_{2\ell},H_{n,k})-N(C_{2\ell},G)
\]
to the minimum number of edge modifications
needed to transform $G$ into a graph isomorphic
to $H_{n,k}$?

\section*{Acknowledgments}
This research was supported by the National Key R\&D Program of China under grant number 2024YFA1013900, the National Natural Science Foundation of China under grant number 12471327.

\section*{Declaration on the use of generative AI}

During the preparation of this work, the authors used ChatGPT (OpenAI)
as an assistive tool to refine proof strategies, check intermediate calculations, and improve the presentation of the manuscript. All mathematical arguments, computations, and conclusions were independently verified and finalized by the authors. The authors
reviewed and edited all AI-assisted content and take full responsibility for the content of this article.

\end{document}